\documentclass[a4paper,9pt]{article}
\usepackage{latexsym,amsfonts,amssymb,amsmath,amscd}
\usepackage{mathrsfs}
\usepackage[amsmath,thmmarks]{ntheorem}
\usepackage{floatflt}
\usepackage[sf]{titlesec}
\usepackage{fancyhdr}
\usepackage{fancyref}
\usepackage{cite}
\usepackage{layouts}
\usepackage{type1cm}
\usepackage{flafter}
\usepackage[perpage,symbol]{footmisc}
\usepackage{indentfirst}
\usepackage{slashbox}
\usepackage{latexsym}

\theoremstyle{plain}
  {\theorembodyfont{\rmfamily}
\theoremstyle{definition}
  {\theorembodyfont{\rmfamily}

\newtheorem{theorem}{Theorem}[section]

\newtheorem{corollary}[theorem]{Corollary}
\newtheorem{proposition}[theorem]{Proposition}

\newtheorem{example}[theorem]{Example}

\newtheorem*{proof}{Proof}
\newtheorem{remark}[theorem]{Remark}

}

\def\deg{\operatorname{deg}}

\def\dgr{\downharpoonright}

\begin{document}

\parindent    = 25pt
\baselineskip = 16pt

\mbox{}\\
\begin{center}
{\Large\bf  Shifted genus expanded  $\cal{W}_{\infty}$ algebra and shifted  Hurwitz numbers}
\vskip 0.4in

Quan Zheng\footnote{E-mail: quanzheng2100@163.com. Partially supported by NSFC.}\qquad

 Mathematics College, Sichuan University\\
610064, Chengdu, Sichuan, PRC

\vskip 0.4in
{\bf Abstract}\\[\baselineskip]

\parbox{14cm}{\quad
 We construct the shifted genus expanded $\cal{W}_{\infty}$ algebra,  which is isomorphic to the central subalgebra $\cal{A}_{\infty}$  of infinite symmetric group algebra and to the shifted Schur symmetrical function algebra $\Lambda^\ast$ defined by A. Y. Okounkov and G. I. Olshanskii. As an application, we get some  differential equations for the generating functions of the shifted
  Hurwitz numbers,
  thus we can express the generating functions in  terms of the shifted genus expanded  cut-and-join operators.
}
\end{center}

\indexspace

\renewcommand\contentsname{\Large Contents}

\tableofcontents

\section{Introduction}
The classical  Hurwitz Enumeration Problem \cite{[H]} and its generalization \cite{[GJ1]},\cite{[GJ2]},\cite{[GJ3]},\cite{[OP]}, etc, have been extensively applied in many mathematical and physical fields, such as the integrability  system theory, the modular space theory,  the relative Gromov-Witten  theory, string theory, referring \cite{[AMMN]}, \cite{[LZZ]}, \cite{[MMN2]}, especially, the investigation of the relationship between the Hurwitz numbers and multi-matrix-model theory have displayed its powerful  vitality, referring \cite{[M1]}, etc. One of the important geometric tools to deal with Hurwitz numbers is the so-called symplectic surgery:
 cutting and gluing \cite{[IP]}, \cite{[L]},\cite{[LR]}, in the views of algebra and differential equations,
 which is equivalent to the so-called cut-and-join operators \cite{[GJ1]}, \cite{[GJ2]}, \cite{[GJ3]},\cite{[LZZ]}, \cite{[OP]}.
 The standard cut-and-join operators can be used to deal with the almost simple Hurwitz numbers
 and the almost simple double Hurwitz numbers  \cite{[GJ1]}, \cite{[GJ2]},\cite{[GJ3]}, \cite{[LZZ]}, \cite{[OP]}.
A. Mironov, A. Morozov, and S. Natanzon, etc, have defined the generalized
 cut-and-join operators in terms of the matrix Miwa variable \cite{[AMMN]}, \cite{[MMN]}, \cite{[MMN1]}.
To distinguish the contributions for the Hurwitz number of
 the source Riemann surface with the different genus, by observing carefully the symplectic surgery (\cite{[L]}, \cite{[LR]}) and the gluing formulas of the
 relative GW-invariants developed by A.M. Li and Y.B. Ruan \cite{[LR]},  and by E. Ionel  and T. Parker  \cite{[IP]}, in \cite{[Z]}, the author has constructed the genus expanded cut-and-join operators and the corresponding genus expanded  differential operator algebra ${\cal{W}}_d$ for any nonnegative integer $d$.
 In this paper, following the idea of \cite{[IK]}, \cite{[MMN2]} and \cite{[OO]}, we extend algebra ${\cal{W}}_d$ to the shifted genus expanded $\cal{W}_{\infty}$ algebra.
 Thus
we have to construct the shifted genus expanded  cut-and-join operators to distinguish the contributions of
the source Riemann surfaces with different genus, i.e., we introduce one parameter $z$
to mark the genus, referring the  formula (\ref{4}). Then after normalizing  the shifted genus expanded  cut-and-join operators by a factor,   we obtain the shifted genus expanded $\cal{W}_{\infty}$ algebra, which is isomorphic to the central subalgebra of infinite symmetric group algebra and to the shifted Schur symmetrical function algebra $\Lambda^\ast$ \cite{[OO]}, referring Corollary \ref{4.6} and Corollary \ref{4.7}.

  As an application of the genus expanded cut-and-join operators, we get some  differential equations for the generating
functions of the shifted Hurwitz numbers  for the source Riemann surface with different genus,
thus we can express the generating functions in  terms of the shifted genus expanded  cut-and-join operators.

\section{  Preliminary: The construction of ${\cal{A}}_{\infty}$ algebra and its structure constants}
First of all, let us introduce some notations.  In the paper, we always arrange that if the variable is negative, then  the any correlated function or invariance    is set to be zero or does NOT appear according to the context. Suppose that $\Delta$ is a Youngian diagram or a partition  of nonnegative integer.
Denote by
$m_i(\Delta)$ the number of rows of length $i$ in the partition
$\Delta$, so we usually denote partition $\Delta$ by $1^{m_1(\Delta)}2^{m_2(\Delta)}3^{m_3(\Delta)}4^{m_4(\Delta)}\cdots$. We denote
$$l(\Delta):=\sum_{i\geq 1}m_{i}(\Delta).$$
$$|\Delta|:=\sum_{i\ge 1}i m_i(\Delta),$$
$$||\Delta||:=\prod_{i\ge 1}i^{m_i(\Delta)},$$
For any two partition $\Delta_1,\Delta_2$, we define their sum and difference as
$$
\Delta_1 \pm \Delta_2=1^{m_1(\Delta_1)\pm m_1(\Delta_2)}2^{m_2(\Delta_1)\pm m_2(\Delta_2)}3^{m_3(\Delta_1)\pm m_3(\Delta_2)}\cdots.
$$
 Let $p=(p_1,p_2,p_3,\cdots,)$ be indeterminantes, which are called the
{\it time-variables}, then we denote
$$p_\Delta:=\prod_{i\geq 1}p_i^{m_i(\Delta)},$$
and
$$\frac{\partial}{\partial p_\Delta}:=\prod_{i\geq 1}\frac{\partial^{m_i(\Delta)}}{\partial p_i^{m_i(\Delta)}},$$

In the section, we recall some basic fact about the construction of ${\cal{A}}_{\infty}$ algebra  and its structure constants following \cite{[IK]}.
Let $\mathbb{P}_n$ be a set of $n$ positive integers $\{1,\ldots,n\}$ , and
 ${\cal{ S}}_n$  be the group of all permutations of  $\mathbb{P}_n$. Fix a subset $E \subset\Bbb P_n$ of size $|E|=r$ and denote by
${\cal{S}}_E$ the group of permutations of the subset $E$.
{\it A partial permutation} of the set $ \mathbb{P}_n$ is a pair
$(E,f)$ consisting of an arbitrary subset
$E\subset\Bbb P_n$ and an arbitrary permutation $f \in {\cal{S}}_E.$
The set $E$ will be referred to as {\it the support} of
$(E,f)$. Define {\it the degree} of a partial permutation $(E,f)\in{\cal{P}}_n$
as $\deg(E,f)=|E|$. Denote by $ {\cal{P}}_n$ the set of all partial permutations of the set
$\Bbb P_n.$
Obviously, the number of elements in ${\cal{P}}_n$ equals
\begin{equation}
 \sum_{k=0}^n {n \choose k}\, k! =
\sum_{k=0}^n (n\dgr k),
\end{equation}
where $n \choose k$ is the binormial coefficient,
and $(n\dgr k)=n(n-1)\ldots(n-k+1)$ is the falling factorial power.

Given two partial permutations $(E_1,f_1)$, $(E_2,f_2)$, we define their
{\it product} as the pair $(E_1\cup
E_2,\,f_1f_2)$, where we can naturally regard $f_1, f_2$ as the permutations of $  {\cal{S}}_{E_1\cup
E_2}$, thus we have the product for $f_1$ and $f_2$. With this multiplication, ${\cal{P}}_n$ becomes a semigroup. The partial
permutation $(\varnothing,e)$, where $e$ is the trivial permutation of
the empty set $\varnothing$, is the unity in ${\cal{P}}_n$.

Denote by ${\cal{B}}_n=\Bbb C[{\cal{P}}_n]$ the complex semigroup algebra of the semigroup
${\cal{P}}_n$, which is semi-simple and  isomorphic to the direct sum of the group
algebras of symmetric groups \cite{[IK]},
$$
 {\cal{B}}_n\cong
\bigoplus_{F\subset\Bbb {\cal{P}}_n} \Bbb C[{\cal{S}}_F].
$$
The centre of the algebra ${\cal{P}}_n$ is of the form
$$
Z({\cal{B}}_n) \cong \bigoplus_{E\subset\Bbb {\cal{P}}_n}
Z(\Bbb C[{\cal{S}}_E]),
$$
where $Z(\Bbb C[{\cal{S}}_E])$ is the centre of the group algebra $\Bbb C[{\cal{S}}_E]$.

The symmetric group ${\cal{S}}_n$ acts on the semigroup
${\cal{P}}_n$ by automorphisms $(E,f)\mapsto(vE,vfv^{-1})$ for $v \in {\cal{S}}_n$.
 The orbits
of this action will be referred to as {\it conjugacy classes} in ${\cal{P}}_n$.
It is obvious that two partial permutations are conjugate if and only if
the sizes of their supports coincide as well as their cycle types. Thus the
conjugacy classes $A_{\Delta;n}\subset{\cal{P}}_n$ are indexed by {\it partial
partitions} of $n$, i.e. by partitions $\Delta\vdash r$ of any
integers $0\le r\le n$. In particular,
$A_{\varnothing;0}=\{(\varnothing,e)\}$. The action of the symmetric group ${\cal{S}}_n$ on ${\cal{P}}_n$ can be continued by
linearity to an action of ${\cal{S}}_n$ on the algebra ${\cal{B}}_n$. Denote by
${\cal{A}}_n={\cal{B}}_n^{{\cal{S}}_n}$ the subalgebra of invariant elements for this action. Let us identify the conjugacy class $A_{\Delta;n}$ with the invariant element
$$
A_{\Delta;n}=\sum_{(E,f)\in A_{\Delta;n}}(E,f)
$$
of the algebra ${\cal{B}}_n$. In particular, if $|\Delta|>n$, then $A_{\Delta;n}=0$.

Given a partial partition $\Delta\vdash r\le n$, denote by
$[\Delta,1^{n-r}]=\Delta\cup\{1^{n-r}\}$, i.e., the partition of $n$ obtained by
adding an appropriate number of unities, which is called a {\it shifting} of $\Delta$. Let
$C_{\Delta}$ be the conjugacy class in the group ${\cal{S}}_{|\Delta|}$ consisting of
permutations of cycle type $\Delta$, and $|C_{\Delta}|$ be the number of the permutations of cycle type $\Delta$.
Let
$C_{\Delta; n}$ be the conjugacy class in the group ${\cal{S}}_n$ consisting of
permutations of cycle type $[\Delta,1^{n-r}]$, and  $|C_{\Delta;n}|$  be the number of the permutations of cycle type $[\Delta,1^{n-r}]$.
Note that $|C_{\Delta}|=|C_{\Delta;n}|$ if $\Delta$ is the partition of $n.$
Denote by
$$\psi:{\cal{P}}_n \rightarrow {\cal{S}}_n$$$$
(E,f)\mapsto f\in {\cal{S}}_E \subseteq {\cal{S}}_n$$
the homomorphism of forgetting
the support of a partial permutation, which obviously can be linearly extended to the  algebra homomorphism  ${\cal{B}}_n \rightarrow {\cal{S}}_n.$
It is clear that
$$
\psi(A_{\Delta;n}) =
{n-|\Delta|+m_1(\Delta) \choose m_1(\Delta)}\, C_{\Delta;n},
$$
thus we can call $\psi(A_{\Delta;n})$ {\it the shifted central element}of ${\cal{S}}_n.$

Assume $m\le n,$ then we introduce a  mapping $\theta_m:{\cal{B}}_n\to{\cal{B}}_m$
by the formula

\begin{equation}
\theta_m(E,f) = \left\{\begin{array}{ccc}
 (E,f)& {\it if} &  E \subset \Bbb{\cal{P}}_m \\
0 &  & otherwise
\end{array}\right.
\end{equation}
The mapping $\theta_m$ is a homomorphism of algebras and it commutes with the
action of the group ${\cal{S}}_m$ on ${\cal{B}}_n$ and on ${\cal{B}}_m$, so we can denote the restricted mapping $\theta_m|_{{\cal{A}}_n}:{\cal{A}}_n\to{\cal{A}}_m$ by the same denotation $\theta_m$.

Denote by ${\cal{B}}_{\infty}$ and ${\cal{A}}_{\infty}$ the projective limit of the algebras ${\cal{B}}_n$ and ${\cal{A}}_n$ with respect to
the morphisms $\theta_n$, respectively.
Let ${\cal{S}}_{\infty}$ be the infinite symmetric group, i.e. the group of finite
permutations of positive integers, thus ${\cal{A}}_{\infty}$ can be regarded as the central subalgebra of group algebra $\mathbb{C}[{\cal{S}}_{\infty}].$

 Denote by $\theta_n$ the natural homomorphism
$\theta_n:{\cal{B}}_{\infty}\to{\cal{B}}_n$ as well as its restriction on ${\cal{A}}_{\infty}$ to ${\cal{A}}_n$.
The natural inclusion of algebras $i_n:{\cal{B}}_n\to{\cal{B}}_{\infty}$ accords with the
projection $\theta_n$: \quad $\theta_n\circ i_n=id_{{\cal{B}}_n}$.

Given a partition $\Delta \vdash r$, let $A_\Delta=\sum
(E,f)$, where the sum extends to partial permutations
$(E,f)\in{\cal{P}}_{\infty}$ such that $|E|=r$ and $f$ has cycle type $\Delta$. The
elements $A_\Delta$, where $\Delta$ runs over all partitions, form a linear basis
in ${\cal{A}}_{\infty}$.
Denote by $\hat{C}_{\Delta_1,\Delta_2}^{\Delta_3}$ the structure constants of the algebra
${\cal{A}}_{\infty}$ in the basis $\{A_\Delta\}$,
\begin{equation}
A_{\Delta_1}\; A_{\Delta_2} =
\sum_{\Delta_3} \hat{C}_{\Delta_1,\Delta_2}^{\Delta_3}A_{\Delta_3}.
\end{equation}
Note that $\theta_n(A_\Delta)=A_{\Delta;n}$, where $A_{\Delta;n}$ is the element
of the algebra ${\cal{A}}_n$. Since
$\theta_n:{\cal{A}}_{\infty}\to{\cal{A}}_n$ and forgetting map $\psi$ are homomorphisms, for any $n\geq 0,$ we have

\proposition \label{prop2.2}(\cite{[IK]}, Proposition 6.1., Theorem 7.1.)
\begin{eqnarray}
A_{\Delta_1;n}\; A_{\Delta_2;n} =
\sum_{\Delta_3}\hat{C}_{\Delta_1,\Delta_2}^{\Delta_3}\, A_{\Delta_3;n},\\
\psi(A_{\Delta_1;n})\;\psi( A_{\Delta_2;n}) =
\sum_{\Delta_3}\hat{C}_{\Delta_1,\Delta_2}^{\Delta_3}\, \psi(A_{\Delta_3;n}).\label{2.6}
\end{eqnarray}

\example \label{2.2} (\cite{[IK]},\cite{[MMN2]})
\begin{eqnarray}
A_{(1)}\,A_{(2)}=2A_{(2)}+A_{(2,1)},\\
A_{(1)}\,A_{(1,1)}=2A_{(1,1)}+3A_{(1,1,1)}\\
A_{(2)}\,A_{(2)}=A_{(1^2)}+3A_{(3)}+2A_{(2^2)}
\end{eqnarray}

\section{Shifted Hurwitz numbers}
Denote  by $\Delta_i=(\delta^i_1,\cdots,\delta^i_{l(\Delta_i)})$, $i=1, \cdots, k,$ a series of  partitions with any degrees.
Let $\Sigma^h$ be a compact (maybe disconnected) Riemann surface of genus $h,$ and $\Sigma^g$
 a compact connected Riemann surface of genus $g$.
 For a given point set $\{q_1,\cdots,q_k\}\in \Sigma^g$, which is called the set of
 branch points, we call a
holomorphic map
$f:\Sigma^h\to \Sigma^g$
a ramified
covering of $\Sigma^g$ of degree $n\ge 0$ by $\Sigma^h$ with a ramification type $(\Delta_{1},\cdots,\Delta_{k})$, if
the preimages of $f^{-1}(q_i)=\{p_i^{1},\cdots,p_i^{l(\Delta_i)},p_i^{l(\Delta_i)+1},\cdots,p_i^{l(\Delta_i)+n-|\Delta_i|}\}$  with orders $ (\delta^i_1,\cdots,\delta^i_{l(\Delta_i)},\underbrace{1,\cdots, 1)}_{n-|\Delta_i|}$ for $i=1,\cdots,k,$ respectively, i.e, we have to shift by many one to the orders of ramification at the branch points  if it is necessary.
Two ramified coverings $f_1$  and $f_2$ with type $(\Delta_{1},\cdots,\Delta_{k})$
 are said to be equivalent if there is a
homeomorphism
$\pi: \Sigma^h\to \Sigma^h$ such that $f_1=f_2\circ\pi$ and $\pi$
preserves the preimages and the ramification type of $f_1$ and $f_2$ at each point $q_i\in \Sigma^g$.
Let $\mu_{g}^{h,n}(\Delta_{1},\cdots,\Delta_{k})$ be the number of equivalent
 covering of $\Sigma^g$ by $\Sigma^h$ with ramification type $([\Delta_{1}, 1^{n-|\Delta_1|}],\cdots,[\Delta_{k}, 1^{n-|\Delta_k|}]),$ which is {\it the classical Hurwitz number}, referring \cite{[LZZ]}.
  We call
 \begin{equation}\label{U}
 U_{g}^{h,n}(\Delta_{1},\cdots,\Delta_{k}):=\prod_{i=1}^{k}{{n-|\Delta_i|+m_1(\Delta_i)} \choose m_1(\Delta_i)}\mu_{g}^{h,n}(\Delta_{1},\cdots,\Delta_{k})
 \end{equation}
  {\it the shifted Hurwitz number}.
  Note that $U_{g}^{h,n}(\Delta_{1},\cdots,\Delta_{k})$ is nonzero only if the Hurwitz formula
\begin{equation}\label{10}
(2-2g)n-(2-2h)=\sum_{i=1}^{k}(|\Delta_i|-l(\Delta_i))
\end{equation}
and inequalities
\begin{equation}
|\Delta_i|\le n, \text{ for } i=1, \cdots, k
\end{equation}
hold. Moreover to stress that the source Riemannian surface is connected, we call  the corresponding shifted Hurwitz number by {\it the connected  shifted Hurwitz number}, which is denoted by $CU_{g}^{h,n}(\Delta_{1},\cdots,\Delta_{k})$.  It is well-known  that the disconnected shifted Hurwitz number can be  related to the connected shifted Hurwitz number by the exponential similar to the classical case, especially, for any given partitions $\Delta_{1},\cdots,\Delta_{k}$ and degree $n$, we have
\begin{eqnarray}\label{CU}
&&U_{g}^{h,n}(\Delta_{1},\cdots,\Delta_{k})p_{\Delta_1}^{(1)}\cdots p_{\Delta_k}^{(k)} \nonumber\\
& \equiv&\sum\prod_i\frac{1}{m_i!}[ CU_{g}^{h_i,n_i}(\Delta_{1}^{(i)},\cdots,\Delta_{k}^{(i)})p_{\Delta_1^i}^{(1)}\cdots p_{\Delta_k^i}^{(k)}]^{m_i},
\end{eqnarray}
where the sum is over all $((n_1,m_1),(n_2,m_2),\cdots)$ and $(\Delta_{1}^{(i)},\cdots,\Delta_{k}^{(i)})$ obeyed the condition:(1) $\sum_j n_jm_j =n;$  (2)$p_{\Delta_1}^{(1)}\cdots p_{\Delta_k}^{(k)}=\prod_i[p_{\Delta_1^i}^{(1)}\cdots p_{\Delta_k^i}^{(k)}]^{m_i},$ and where $p_{\Delta_1}^{(1)}, \cdots,  p_{\Delta_k}^{(k)}$, $p_{\Delta_1^i}^{(1)}\cdots p_{\Delta_k^i}^{(k)}$ are the time-variables, moreover $h_i$ satisfies the following equality:
\begin{equation}
(2-2g)n_i-(2-2h_i)=\sum_{j=1}^{k}(|\Delta_j^{(i)}|-l(\Delta_j^{(i)})).
\end{equation}
By \cite{[H]}, using the notations of section 2, we have
\begin{eqnarray}
U_{g}^{h,n}(\Delta_{1},\cdots,\Delta_{k})=\frac{1}{n!}[1]
\prod_{j=1}^g\prod_{a_j,b_j\in {\cal{S}}_n}[a_j,b_j]\psi({A}_{\Delta_1;n})\cdots \psi({A}_{\Delta_k;n}),
\end{eqnarray}
i.e., $U_{g}^{h,n}(\Delta_{1},\cdots,\Delta_{k})$ equals to $\frac{1}{n!}$ times the
coefficient of the identity of the product of $g$-tuple  commutators $\prod_{a,b\in S_n}[a,b]$ and $k$-tuple
  shifted central elements $\psi(A_{\Delta_1;n}), \cdots, \psi(A_{\Delta_k;n})$,
    thus, according to the theory of representations of the symmetric group ${\cal{S}}_n$, which is also equivalent to (referring  \cite{[AMMN]})
\begin{eqnarray}{\label{15}}
U_{g}^{h,n}(\Delta_{1},\cdots,\Delta_{k})=
\sum_{\lambda \vdash n}({\frac {dim{\lambda}}{n!}})^{2-2g}\phi_{\lambda}(\Delta_1)\cdots\phi_{\lambda}(\Delta_k),
\end{eqnarray}
which expresses them through the properly normalized  characters
\begin{equation}\label{21}
\phi_\lambda(\Delta) \equiv \left\{\begin{array}{ccc}
0 & {\rm for} &  |\Delta|>|\lambda| \\
{|\lambda|-|\Delta|+m_1(\Delta)\choose m_1(\Delta)}\frac{1}{dim\lambda}|C_{\Delta; |\lambda|}|\chi_\lambda(\Delta,1^{|\lambda|-|\Delta|})
 & {\rm for} &  |\Delta|\leq |\lambda|
\end{array}\right.
\end{equation}
where   $\lambda$ is a Young diagram of degree $|\lambda|$, and $dim\lambda$ is the dimension of the irreducible  representation of the symmetric group ${\cal{S}}_{|\lambda|}$ corresponding to $\lambda$, $
\chi_\lambda(\Delta,1^{|\lambda|-|\Delta|})
$  is the  character of a permutation $f\in C_{\Delta,|\lambda|} \subset \mathbb{C}[{\cal{S}}_{|\lambda|}]$ under the irreducible  representation $\lambda$. Let $p=(p_1,p_2,p_3,\cdots,)$ be  time-variables,
then for any partition $\lambda,$ it is well known that the normalized
characters $\phi_{\lambda}(\Delta)$ are related to the Schur functions $S_{\lambda}\{p\}$
as \cite{[M]}, \cite{[MMN]}:
\begin{equation}\label{5}
S_{\lambda}\{p\} = \sum_{\Gamma^{\prime}}
\frac{ dim{\lambda}}{|\lambda|!}\phi_{\lambda}(\Gamma^{\prime})p_{\Gamma^{\prime}}\delta_{|\lambda|,|\Gamma^{\prime}|},
\end{equation}
or to the shifted Schur functions  $S_{\lambda}\{p_m+\delta_{m,1}\}$ as \cite{[MMN]}:
\begin{equation}
S_{\lambda}\{p_m+\delta_{m,1}\} = \sum_{\Gamma^{\prime}}
\frac{ dim{\lambda}}{|\lambda|!}\phi_{\lambda}(\Gamma^{\prime})p_{\Gamma^{\prime}},
\end{equation}
where $S_{\lambda}\{p_m+\delta_{m,1}\}$ means that $S_{\lambda}$ is the function of $p_1+1, p_2, p_3, \cdots.$

The following examples is important to determine the lower degree genus expanded cut-and-join operators, see Example \ref{26}:
\example For any positive integer $a, b,$ we have
\begin{eqnarray}
CU_{0}^{0,1}((0),(1),(0))=1;\label{17}\\
 CU_{0}^{0,a}((a),(0),(a))=\frac{1}{a}\\
 CU_{0}^{0,a}((a),(1),(a))=1;\\
  CU_{0}^{0,a}((a),(1,1),(a))=\frac{1}{2}(a-1)\\
   CU_{0}^{0,a}((a),(1,1,1),(a))=\frac{1}{6}(a-1)(a-2)\\
  CU_{0}^{0,a+b}((a,b),(2),(a+b))=1-\frac{\delta_{a,b}}{2}\\
  CU_{0}^{0,a+b}((a,b),(2,1),(a+b))=(1-\frac{\delta_{a,b}}{2})(a+b-2)
\end{eqnarray}
\bigskip
\remark From the above examples, we note that the relative Gromov-Witten invariants even with tangent multiple one is different to the relative Gromov-Witten invariants without tangent conditions.

\section{Shifted genus expanded   ${\cal{W}}_{\infty}$ algebra}
For every nonnegative integer $n,$ in \cite{[Z]}, we constructed
a genus expanded differential algebra ${\cal{W}}_n$, which is algebraic isomorphic to the central  subalgebra $Z(\mathbb{C}[{\cal{S}}_n])$ (referring  \cite{[Z]}, Corollary 3.6). In this section we will construct shifted genus expanded differential algebra ${\cal{W}}_{\infty}$ by observing the symplectic surgery and the construction of ${\cal{A}}_\infty$ algebra in  the  section 2.

 For any partition $\Delta$,  we call a series of partitions $\tilde{\Delta}=(\Delta^1,\Delta^2,\cdots)$ {\it a proper Re-partition} of $\Delta$ if $\Delta= \sum_{i}\Delta^i$ and $|\Delta^i| \geq 1$, for any $i.$  We denote by ${\cal{PP}}_{\Delta}$ the all proper Re-partitions of $\Delta.$  For any series partition   $(\Gamma_1^{\prime}, \Gamma_2^{\prime}, \cdots)$ and $(\Gamma_1, \Gamma_2, \cdots)$, denote  $\sum_i\Gamma^{\prime}_i:=\Gamma^{\prime}$ and $ \sum_i\Gamma_i:=\Gamma$. We define $|Aut(\Gamma^{\prime}, \tilde{\Delta}, \Gamma)|$ the number of the automorphisms of the triples: $(\Gamma_i^{{\prime}},\Delta^{i},\Gamma_i), i=1, 2, \cdots$.

 At now,  we  can associate any partition $\Delta$
 a {\it shifed genus  expanded cut-and-join differential  operator} $W(\Delta,z)$ as follows:
\begin{eqnarray}\label{4}
&&W(\Delta,z)=\\
&\sum_{\tilde{\Delta}\in {\cal{PP}}_{\Delta}}&\sum_{\Gamma^{\prime},\Gamma}\frac{z^{|\Delta|-l(\Delta)+l(\Gamma^{\prime})
-l(\Gamma)}}
{|Aut(\Gamma^{\prime},\tilde{\Delta}, \Gamma)| }:\prod_{i}[||\Gamma_i^{\prime}||\delta_{|\Gamma_i|,|\Gamma_i^{\prime}|}CU_{0}^{h_i^+,|\Gamma_i^{\prime}|}(\Gamma_i^{\prime},\Delta^i,\Gamma_i)
p_{\Gamma_i}\frac{\partial}{\partial p_{\Gamma_i^{\prime}}}]:,\nonumber
\end{eqnarray}
where : : means that the classical normal ordering product, i.e, all $p_i$  will always  be set in the front of the all partial operator $\frac{\partial}{\partial p_j},$
and  genus $h^i_+ \geq 0$ is determined  by the following:
\begin{equation}
2h_i^+-2=|\Delta^i|-l(\Delta^i)-l(\Gamma_i^{\prime})-l(\Gamma_i).
\end{equation}
\bigskip

\remark
In generally, the following equality:
\begin{equation}\label{25}
U_{0}^{h_+,|\Gamma^{\prime}|}(\Gamma^{\prime},\Delta,\Gamma)=\sum_{\tilde{\Delta}}\frac{1}{|Aut(\Gamma^{\prime},\tilde{ \Delta}, \Gamma)| }\prod_{i}[\delta_{|\Gamma_i^{\prime}|,|\Gamma_i|}CU_{0}^{h_i^+,|\Gamma_i^{\prime}|}(\Gamma_i^{\prime},\Delta^i,\Gamma_i)]
\end{equation}
doesn't hold due to   shifting the ramification type,  for example: $$U_{0}^{-1,7}((4,3),(2,1),(4,2,1))=\frac{5}{4},$$
 but\\ $$CU_{0}^{0,4}((4),(1),(4))CU_{0}^{0,3}((3),(2),(2,1))=1,$$ referring  (\cite{[MMN]}, Section 2.2).  We denote by $U(\Gamma^{\prime},\Delta,\Gamma)$ the RHS of the formula (\ref{25}).
\bigskip

\bigskip
\example \label{26}
\begin{eqnarray}
W((1),z)=\sum_aap_a\frac{\partial}{\partial p_a},\\
W((1,1),z)=\frac{1}{2}\sum_a a(a-1)p_a\frac{\partial}{\partial p_{a}}+\frac{1}{2}\sum_{a,b}abp_ap_b\frac{\partial^2}{\partial p_{a}\partial p_b},\\
W((1,1,1))=
{1\over 6}\sum_{a} a(a-1)(a-2)p_a\frac{\partial}{\partial p_a}\,\nonumber\\
+{1\over 2}\sum_{a,b} a(a-1)bp_ap_b\frac{\partial^2}{\partial p_a\partial p_b}\,
+{1\over 6}\sum_{a,b,c}abcp_ap_bp_c{\partial^3\over\partial p_a\partial p_b\partial p_c},\\
W((2),z)=\frac{1}{2}\sum_{a,b}(a+b)p_ap_b\frac{\partial}{\partial p_{a+b}}+\frac{1}{2}\sum_{a,b}z^{2}abp_{a+b}\frac{\partial^2}{\partial p_{a}\partial p_b}, \label{22}\\
W((2,1),z)=
\frac{1}{2}\sum_{a,b}(a+b)(a+b-2)p_ap_b\frac{\partial}{\partial p_{a+b}}+\frac{1}{2}\sum_{a,b,c}(a+b)cp_ap_bp_c\frac{\partial^2}{\partial p_{a+b} \partial p_c}\nonumber\\\label{23}
+\frac{1}{2}\sum_{a,b,c}z^2abcp_ap_{b+c}\frac{\partial^3}{\partial p_a \partial p_b\partial p_c}+
\frac{1}{2}\sum_{a,b}z^{2}ab(a+b-2)p_{a+b}\frac{\partial^2}{\partial p_a \partial p_b}.
\end{eqnarray}

\remark The formula (\ref{22}) is the standard cut-and-join operator, referring \cite{[GJ1]},\cite{[GJ2]}, \cite{[GJ3]}, \cite{[LZZ]}, \cite{[OP]}.

\bigskip
To give out the eigenfunctions of the universal genus expanded cut-and-join operator  $W(\Delta,z)$,
we define {\it the ``genus expanded"
Schur functions} $S_{\lambda}\{p,z\}$ similar to formula (\ref{5})
as follows:
\begin{equation}\label{11}
S_{\lambda}\{p,z\} := \sum_{\Gamma^{\prime}} z^{-|\Gamma^{\prime}|-l(\Gamma^{\prime})}\frac{ dim{\lambda}}{|\lambda|!}\phi_{\lambda}(\Gamma^{\prime})p_{\Gamma^{\prime}}\delta_{|\lambda|,|\Gamma^{\prime}|}.
\end{equation}

\bigskip
\theorem \label{4.4}
  For any partitions $\Delta_1, \Delta_2,$ as operators on the functions of the time-variables
$p = (p_1, p_2,\cdots )$,  we have
\begin{equation}\label{19}
{W(\Delta_1,z) W(\Delta_2,z) = \sum_{\Delta_3}z^{(|\Delta_1|-l(\Delta_1))+(|\Delta_2|-l(\Delta_2))+(-|\Delta_3|+l(\Delta_3)}
\hat{C}_{\Delta_1\Delta_2}^{\Delta_3}  W(\Delta_3,z) },
\end{equation}\label{CWW0}
where
$\hat{C}_{\Delta_1\Delta_2}^{\Delta_3}$
is the  structure constants
of  ${\cal{A}}_\infty$. Moreover, we have
\begin{equation}\label{27}
 W(\Delta,z)S_\lambda\{p,z\} = z^{|\Delta|-l(\Delta)}\phi_\lambda(\Delta) S_\lambda\{p,z\}.
\end{equation}
\proof
 Let us firstly prove the equality (\ref{27}).
For any partition $\Delta^{\prime},$ we have
\begin{eqnarray}
&&W(\Delta,z)p_{\Delta^{\prime}}=\\
&&\sum_{\Gamma^{\prime},\Gamma}z^{|\Delta|-l(\Delta)+l(\Gamma^{\prime})
-l(\Gamma)}\delta_{|\Gamma|,|\Gamma^{\prime}|}||\Gamma^{\prime}||
U(\Gamma^{\prime},\Delta,\Gamma)
p_{\Gamma}\prod_{i\geq1}(m_i(\Delta^{\prime})\dgr {m_i(\Gamma^{\prime})}p_i^{m_i(\Delta^{\prime})-m_i(\Gamma^{\prime})}.\nonumber
\end{eqnarray}
Assume $\Delta^{\prime\prime}:=\Gamma+(\Delta^{\prime}-\Gamma^{\prime})$, then we have
\begin{eqnarray}\label{39}
&&W(\Delta,z)p_{\Delta^{\prime}}\nonumber\\
&=&\sum_{\Delta^{\prime\prime}}z^{|\Delta|-l(\Delta)+l(\Delta^{\prime})
-l(\Delta^{\prime\prime})}p_{\Delta^{\prime\prime}}\sum_{\Gamma^{\prime}}\delta_{|\Gamma|,|\Gamma^{\prime}|}||\Gamma^{\prime}||
U(\Gamma^{\prime},\Delta,\Delta^{\prime\prime}-(\Delta^{\prime}-\Gamma^{\prime}))
\prod_{i\geq1}(m_i(\Delta^{\prime})\dgr {m_i(\Gamma^{\prime})}\nonumber\\
&=&\sum_{\Delta^{\prime\prime}}z^{|\Delta|-l(\Delta)+l(\Delta^{\prime})
-l(\Delta^{\prime\prime})}p_{\Delta^{\prime\prime}}||\Delta^{\prime}||
U_0^{h^+,|\Delta^{\prime\prime}|}(\Delta^{\prime},\Delta,\Delta^{\prime\prime})\delta_{|\Delta^{\prime}|,|\Delta^{\prime\prime}|},
\end{eqnarray}
where the  second  equality follows from the fact that the shifted Hurwitz number $U_0^{h^+,|\Delta^{\prime\prime}|}(\Delta^{\prime},\Delta,\Delta^{\prime\prime})$ can be expressed in the many connected shifted Hurwitz number which especially include the form  $CU_0^{0,a}(a,(0),a)$, which corresponds to the divisor $p_a$ that need not to be derived, referring the formula (\ref{CU}). For example, let $\Delta=(2,1)$, $\Delta^{\prime}=(4,3)$, $\Delta^{\prime\prime}=(4,2,1),$ then we have
\begin{eqnarray}
&&12U_{0}^{-1,7}((4,3),(2,1),(4,2,1))\nonumber\\
&=&12CU_{0}^{0,4}((4),(1),(4))CU_{0}^{0,3}((3),(2),(2,1))\nonumber\\
&&+12CU_{0}^{0,4}((4),(0),(4))CU_{0}^{0,3}((3),(2,1),(2,1)).
\end{eqnarray}
Thus the equality (\ref{27})  follows from the formula (\ref{11}), (\ref{39}), (\ref{15}) and the orthogonal  relation of the irreducible characteristic  of symmetric group (\cite{[Z]}, Lemma 2.1) as following:
\begin{eqnarray}
&&W(\Delta,z)S_{\lambda}\{p,z\} \nonumber\\
&=& \sum_{\Delta^{\prime}} \frac{ dim{\lambda}}{|\lambda|!}\phi_{\lambda}(\Delta^{\prime})\delta_{|\Delta^{\prime}|,|\lambda|}\sum_{\Delta^{\prime\prime}}z^{|\Delta|-l(\Delta)-|\Delta^{\prime\prime}|-l(\Delta^{\prime\prime})}p_{\Delta^{\prime\prime}}||\Delta^{\prime}||
U_0^{h^+,|\Delta^{\prime\prime}|}(\Delta^{\prime},\Delta,\Delta^{\prime\prime})\delta_{|\Delta^{\prime}|,|\Delta^{\prime\prime}|}\nonumber\\
&=&\sum_{\Delta^{\prime\prime}} z^{|\Delta|-l(\Delta)-|\Delta^{\prime\prime}|-l(\Delta^{\prime\prime})}\sum_{\Delta^{\prime}}\frac{ dim{\lambda}}{|\lambda|!}\phi_{\lambda}(\Delta^{\prime})\delta_{|\Delta^{\prime}|,|\lambda|}||\Delta^{\prime}||\times\nonumber\\
&&\sum_{|\mu|=|\Delta^{\prime\prime}|}(\frac{ dim{\mu}}{|\mu|!})^2\phi_{\mu}(\Delta^{\prime})\phi_{\mu}(\Delta)\phi_{\mu}(\Delta^{\prime\prime})p_{\Delta^{\prime\prime}}\delta_{|\Delta^{\prime}|,|\Delta^{\prime\prime}|}\nonumber\\
&=&\sum_{\Delta^{\prime\prime}} z^{|\Delta|-l(\Delta)-|\Delta^{\prime\prime}|-l(\Delta^{\prime\prime})}\sum_{\mu}\frac{ dim{\lambda}}{|\lambda|!}\delta_{\mu,\lambda}\phi_{\mu}(\Delta)\phi_{\mu}(\Delta^{\prime\prime})p_{\Delta^{\prime\prime}}\delta_{|\Delta^{\prime\prime}|,|\mu|}\nonumber\\
&=& z^{|\Delta|-l(\Delta)}\phi_\lambda(\Delta) S_\lambda\{p,z\}.
\end{eqnarray}
The formula (\ref{19}) follows from the following equalities:
\begin{eqnarray}
&&W(\Delta_1,z) W(\Delta_2,z)S_\lambda\{p,z\} \nonumber\\
&=& z^{|\Delta_2|-l(\Delta_2)}\phi_\lambda(\Delta_2) W(\Delta_1,z)S_\lambda\{p,z\}\nonumber\\
&=&z^{|\Delta_1|-l(\Delta_1)+|\Delta_2|-l(\Delta_2)}\phi_\lambda(\Delta_1)\phi_\lambda(\Delta_2)S_\lambda\{p,z\}\nonumber\\
&=&z^{|\Delta_1|-l(\Delta_1)+|\Delta_2|-l(\Delta_2)}\hat{C}_{\Delta_1\Delta_2}^{\Delta_3}  \phi(\Delta_3)S_\lambda\{p,z\}\\
&=&z^{|\Delta_1|-l(\Delta_1)+|\Delta_2|-l(\Delta_2)-(|\Delta_3|-l(\Delta_3))}\hat{C}_{\Delta_1\Delta_2}^{\Delta_3}  W(\Delta_3,z)S_\lambda\{p,z\},\nonumber
\end{eqnarray}
where the third equality   comes from  (\cite{[MMN]}, formula (45)).

\bigskip

\example By Example \ref{26}, it is easy to check that
\begin{eqnarray}
W((1),z)W((2),z)=2W((2),z)+W((2,1),z),\\
W((1),z)W((1,1),z)=2W((1,1),z)+3W((1,1,1),z),
\end{eqnarray}
which coincide with the Example \ref{2.2} and Theorem \ref{4.4}.
\corollary \label{4.6}
If we normalize the shifted genus expanded cut-and-join operator $W(\Delta, z)$ by a factor  $z^{-|\Delta|+l(\Delta)}$
\begin{equation}\label{20}
\hat{W}(\Delta, z):=z^{-|\Delta|+l(\Delta)}W(\Delta, z),
\end{equation}
then  as operators on the space of functions in time-variables $p=(p_1, p_2, \cdots)$,
all shifted genus expanded  cut-and-join operators $\hat{W}(\Delta, z)$  form a commutative associative
algebra, denoted by ${\cal{W}}_\infty,$ which is called by {\it the shifted genus expanded  algebra. }
 \begin{equation}
 {\hat{W}(\Delta_1,z) \hat{W}(\Delta_2,z) = \sum_{\Delta_3}
\hat{C}_{\Delta_1\Delta_2}^{\Delta_3 } \hat{W}(\Delta_3,z) }
\end{equation}\label{CWW0}
 i.e., we have an algebraical isomorphism:
\begin{eqnarray}
{\cal{W}}_\infty &\cong &{\cal{A}}_\infty\nonumber\\
        \hat{W}(\Delta,z)   & \mapsto & A_{\Delta}.
\end{eqnarray}
Moreover,  $\hat{W}(\Delta,z)$ have the genus expanded  Schur function  $S_\lambda\{p,z\}$
 as their eigenfunctions
and $\phi_\lambda(\Delta)$ as the corresponding eigenvalues:
\begin{equation}
 \hat{W}(\Delta,z)S_\lambda\{p,z\} = \phi_\lambda(\Delta) S_\lambda\{p,z\}.
\end{equation}
\proof By straightforward  calculation, we omit it.

\corollary \label{4.7} The shifted genus expanded $\cal{W}_{\infty}$ algebra is isomorphic to the shifted Schur symmetrical function algebra $\Lambda^\ast$  defined by A. Y. Okounkov, and G. I. Olshanskii \cite{[OO]}.
\proof By (\cite{[IK]}, Theorem 9.1), we know that  ${\cal{A}}_\infty$ is isomorphic to  the shifted Schur symmetrical function algebra $\Lambda^\ast$, thus the result follows from Corollary \ref{4.6}.

\remark Due to the Corollary \ref{4.7}, combined the result of the section 3, we  call the genus expanded cut-and-join operator $W(\Delta,z)$ and its normalization $\hat{W}(\Delta,z)$ as the {\it shifted genus expanded cut-and-join operators}.
\remark Note a fact that the  structure constants of the shifted genus expanded  algebra  ${\cal{W}}_\infty$ do not dependent on the parameter $z$ although the shifted genus expanded cut-and-join operators do, which is because that $z$  recorded the $``lost" $ genus when we  execute symplectic surgery and  we can obtain the same results  if we execute the symplectic surgery once instead of twice, referring \cite{[Z]}, Remark 4.3.

\section {Generating functions of the shifted Hurwitz number and its differential equations}
To apply the shifted genus expanded cut-and-join operators, we assme that  $\Delta_1, \cdots,\Delta_m$ are any fixed partitions,  then we can define a
generating function
\begin{eqnarray}\label{5.1}
&&\Phi_g\{z|(u_1,\Delta_1),\cdots,(u_m,\Delta_m)|p^{(1)},\cdots,p^{(k)},p\} \nonumber\\
&=&\sum_{n\geq 0} \sum_{l_1,\cdots,l_m\ge 0}
\sum_{\Gamma,\Gamma_1,\ldots,\Gamma_{k}}z^{2h-2}U_{g}^{h,n}(\underbrace{\Delta_{1},\cdots,\Delta_1}_{l_1},\cdots,\underbrace{\Delta_{m},\cdots,\Delta_m}_{l_m},\Gamma_1,\cdots,
\Gamma_k,\Gamma)[\prod_{j=1}^{m}\frac{u_j^{l_j}}{l_j!}][\prod_{i=1}^{k}p^{(i)}_{\Gamma_i}]p_\Gamma \nonumber\\
&=&\sum_{n\geq 0}\sum_{|\lambda|=n}\sum_{l_1,\cdots,l_m\ge 0}\sum_{\Gamma,\Gamma_1,\ldots,\Gamma_{k}}z^{2h-2}({\frac {dim{\lambda}}{|\lambda|!}})^{2-2g}[\prod_{\makebox {j=1}}^m(\phi_{\lambda}(\Delta_j))^{l_j}\frac{(u_j)^{l_j}}{l_j!}]
[\prod_{i=1}^{k}\phi_{\lambda}(\Gamma_i)p^{(i)}_{\Gamma_i}]\phi_{\lambda}(\Gamma)p_\Gamma,
\end{eqnarray}
where $z, u_1,\cdots, u_m$ are indeterminate variables, $p, p^{(1)}, \cdots, p^{(k)}$ are time-variables, and
$2h-2$ is determined  by the Hurwitz formula:
\begin{equation}\label{g1}
(2-2g)n-(2-2h)=\sum_{j=1}^{m}l_j(|\Delta_j|-l(\Delta_j))+\sum_{j=1}^{k}(|\Gamma_j|-l(\Gamma_j))+(|\Gamma|-l(\Gamma)).
\end{equation}
Moreover, we have some special initial values \cite{[M]}, \cite{[MMN]}:
\item
\begin{eqnarray}\label{35}
\Phi_0\{z||p \}&=&\sum_{\lambda}\sum_{\Delta}z^{-|\Delta|-l(\Delta)}(\frac{dim{\lambda}}{|\lambda|!})^2\phi_{\lambda}(\Delta)p_{\Delta}\nonumber\\
&=&\sum_{\lambda}\frac{dim{\lambda}}{|\lambda|!}S_{\lambda}\{p_m+\delta_{m,1},z\}\nonumber\\
&=&exp{\frac{p_1+1}{z^2}};
\end{eqnarray}
\begin{eqnarray}\label{2}
\Phi_0\{z||p^{(1)},p\}&=&\sum_{\lambda}\sum_{\Delta_1,\Delta_2}z^{2|\lambda|-|\Delta_1|-|\Delta_2|-l(\Delta_1)-l(\Delta_2)}(\frac{dim{\lambda}}{d!})^2\phi_{\lambda}(\Delta_1)\phi_{\lambda}(\Delta_2)p^{(1)}_{\Delta_1}p_{\Delta_2}\nonumber\\
&=&\sum_{\lambda}z^{2|\lambda|}S_{\lambda}\{p^{(1)}_m+\delta_{m,1},z\}S_{\lambda}\{p_m+\delta_{m,1},z\}\nonumber\\
&=&exp(\sum_{m\geq 1}\frac{ (p^{(1)}_m+\delta_{m,1})(p_m+\delta_{m,1})}{z^2})
\end{eqnarray}
where we have to add one to $p^{(1)}_1$ and $p_1$ due to the shifting, referring the equality (\ref{17}).

\remark
In the formula (\ref{5.1}), we sum all the degree of the covering map $n \geq 0$, which is different from \cite{[Z]} formula (23).
\bigskip
\theorem For any $i$, we have
\begin{eqnarray}\label{3}
&&\nonumber \frac{\partial\Phi_g\{z|(u_1,\Delta_1),\cdots,(u_m,\Delta_m)|p^{(1)},\cdots,p^{(k)},p\}}{\partial u_i} \nonumber\\
& =&W(\Delta_i,z)\Phi_g\{z|(u_1,\Delta_1),\cdots,(u_m,\Delta_m)|p^{(1)},\cdots,p^{(k)},p\}.
\end{eqnarray}
\proof
Obviously, we have
$$
\frac{\partial\Phi_g\{z|(u_1,\Delta_1),\cdots,(u_n,\Delta_n)|p^{(1)},\cdots,p^{(k)},p\}}{\partial u_i}$$
$$=\sum_{n\geq 0}\sum_{l_1,\cdots,l_n\ge 0}
\sum_{\Gamma_1,\ldots,\Gamma_{k}, \Gamma}z^{2h-2}\mu_{g}^{h,n}(\underbrace{\Delta_{1},\cdots,\Delta_1}_{l_1},\cdots,\underbrace{\Delta_{m},\cdots,\Delta_m}_{l_m},\Gamma_1,\cdots,
\Gamma_k,\Gamma)$$
$$\frac{(u_i)^{l_i-1}}{(l_i-1)!}[\prod_{\makebox {j=1,}\\ j\ne i}^m\frac{(u_j)^{l_j}}{l_j!}]
[\prod_{j=1}^{k}p^{(j)}_{\Gamma_j}]p_\Gamma.
$$
We can write RHS of equation (\ref{3}) as
\begin{eqnarray*}
\mathrm{RHS}&=&\sum_{n\geq 0}\sum_{l_i\ge 1} \sum_{l_1,\cdots,\check{l_i},\cdots,l_m\ge 0}
\sum_{\Gamma_1,\ldots,\Gamma_{k},\Gamma^{\prime}}z^{2h^i_--2}\\
&&\times
\mu_{g}^{h^i_-,n}(\underbrace{\Delta_{1},\cdots,\Delta_1}_{l_1},\cdots,\underbrace{\Delta_{i},
\cdots,\Delta_i}_{l_i-1},\cdots,\underbrace{\Delta_{m},\cdots,\Delta_m}_{l_m},\Gamma_1,\cdots,
\Gamma_k,\Gamma^{\prime})\\
&&\times
\frac{(u_i)^{l_i-1}}{(l_i-1)!}[\prod_{\makebox {j=1,} j\ne i}^n\frac{(u_j)^{l_j}}{l_j!}]
[\prod_{j=1}^{k}p^{(j)}_{\Gamma_j}]W(\Delta_i,z)p_{\Gamma^{\prime}},
\end{eqnarray*}
where $\check{l_i}$ means that we omit $l_i$, and  $2h^i_--2$ is also determined by the Hurwitz formula:
\begin{equation}\label{g2}
(2-2g)n-(2-2h^i_-)=\sum_{j=1}^{m}l_j(|\Delta_j|-l(\Delta_j))-(|\Delta_i|-l(\Delta_i))+\sum_{j=1}^{k}(|\Gamma_j|-l(\Gamma_j))+(|\Gamma^{\prime}|-l(\Gamma^{\prime})),
\end{equation}

Then the theorem follows from the following facts:
\begin{itemize}
\item Fact (1): (referring the formula (\ref{39}))
\begin{eqnarray}
&& W(\Delta,z)p_{\Delta^{\prime}}\\
&=&\sum_{\Delta^{\prime\prime}}z^{|\Delta|-l(\Delta)+l(\Delta^{\prime})
-l(\Delta^{\prime\prime})}p_{\Delta^{\prime\prime}}||\Delta^{\prime}||
U_0^{h^+,|\Delta^{\prime\prime}|}(\Delta^{\prime},\Delta,\Delta^{\prime\prime})\delta_{|\Delta^{\prime}|,|\Delta^{\prime\prime}|}.\nonumber
\end{eqnarray}
\item  Fact (2):
 $CU_{0}^{h^i_+,|\beta|}(\alpha,\Delta^i,\beta)\delta_{|\alpha|,|\beta|}\ne 0$ only if
\begin{equation}\label{g3}
2h^i_+-2=-l(\alpha)+(|\Delta^i|-l(\Delta^i))-l(\beta).
\end{equation}
\end{itemize}

 \bigskip
Immediately, we have
\corollary
\begin{equation}\label{38}
\Phi_g\{z|(u_1,\Delta_1)\cdots,(u_n,\Delta_n)|p^{(1)},\cdots,p^{(k)},p\} =[\prod_{i=1}^n exp (u_iW(\Delta_i,z))]
\Phi_g\{z||p^{(1)},\cdots,p^{(k)},p\}
\end{equation}

\example \label{5.3}
Assume  $\Delta=(2,1),$ then  by the formula (\ref{23}), (\ref{35}) and (\ref{38}),
 we get the  generating function(with any genus) :
\begin{eqnarray}
&&\Phi_0\{z|(u,(2,1))|p\}=\\
&&\frac{1}{6}(p_1+1)^3z^{-6}+ \frac{1}{2}u(p_1+1)p_2z^{-4}+\frac{1}{2}u^2p_3z^{-2}+\frac{1}{4}u^2(p_1+1)^3z^{-4}+\frac{3}{4}u^3(p_1+1)p_2z^{-2}+\cdots,\nonumber
\end{eqnarray}
where we omit the higher terms of $u$ and $n.$

\bigskip
{\bf Acknowledgements} The author  would like to thank the every members of geometry
team of Mathematics College in Sichuan University, and Prof. Yongbin Ruan, Prof. Yuping Tu,
Prof. Guohui Zhao, Prof. Qi Zhang for their helpful discussions.

\end{document}